\documentclass[12pt]{article}
\usepackage{amsmath,amsfonts,latexsym,amsthm,amssymb}
\newcommand{\D}{\mathbb D_t^{(\alpha )}}
\newcommand{\Rn}{\mathbb R^n}
\DeclareMathOperator{\const}{const}
\DeclareMathOperator{\R}{Re}

\begin{document}
\newtheorem{prop}{Proposition}
\newtheorem*{lem}{Lemma}
\newtheorem{teo}{Theorem}
\pagestyle{plain}
\title{Asymptotic Properties of Solutions of the Fractional Diffusion-Wave Equation}
\author{Anatoly N. Kochubei
\\ \footnotesize Institute of Mathematics,\\
\footnotesize National Academy of Sciences of Ukraine,\\
\footnotesize Tereshchenkivska 3, Kiev, 01601 Ukraine\\
\footnotesize E-mail: \ kochubei@i.com.ua}
\date{}
\maketitle

\vspace*{3cm}
\begin{abstract}
For the fractional diffusion-wave equation with the Caputo-Dzhrbashyan fractional derivative of order $\alpha \in (1,2)$ with respect to the time variable, we prove an analog of the principle of limiting amplitude (well-known for the wave equation and some other hyperbolic equations) and a pointwise stabilization property of solutions (similar to a well-known property of the heat equation and some other parabolic equations).
\end{abstract}
\vspace{2cm}


\newpage
\section{Introduction}

The fractional diffusion-wave equation
\begin{equation}
\left( \D u\right) (t,x)-\Delta u(t,x)=f(t,x),\quad t\in (0,T],\ x\in
\mathbb R^n,
\end{equation}
where $1<\alpha <2$, $\D$ is the Caputo-Dzhrbashyan fractional derivative, that is
$$
\left( \mathbb D_t^{(\alpha )}u\right) (t,x)=\frac{1}{\Gamma (2-\alpha )}\frac{\partial}{\partial t}\int\limits_0^t (t-\tau )^{-\alpha +1}u'_\tau (\tau ,x)\,d\tau -t^{-\alpha +1}\frac{u'_t(0,x)}{\Gamma (2-\alpha )},
$$
describes the propagation of stress pulses in a viscoelastic medium \cite{Ma}. Its properties are intermediate between those of the classical heat and wave equations. This equation and its generalizations have been studied by many authors (see \cite{K13,K14,Ma,Ps} for further references) who put their emphasis either on properties similar to those of parabolic equations, like the regularity properties, or those resembling hyperbolic equations, like the exponential decay of the fundamental solution outside the fractional light cone.

Note that we use an expression for the fractional derivative of order from the interval $(1,2)$, which is equivalent to the usual one (see the proof of Theorem 2.1 in \cite{KST}) and contains explicitly the first derivative whose continuity is assumed for classical solutions of the equation (1).

In this paper devoted to asymptotic properties of solutions of the equation (1), we study two kinds of asymptotic behavior typical for hyperbolic and parabolic equations, respectively -- the principle of limiting amplitude and the stabilization property. We find that the fractional diffusion-wave equation possesses both of them simultaneously.

The principle of limiting amplitude is a classical property of the wave equation on $(0,\infty )\times \mathbb R^3$ \cite{TS} extended later to more general equations and domains: if $u(t,x)$ is a solution of the equation $\dfrac{\partial^2 u}{\partial t^2}-\Delta u=e^{-i\mu t}f(x)$ ($f$ is a function with compact support) with zero initial conditions, then $u(t,x)=e^{-i\mu t}v(x)+o(1)$, $t\to \infty$, where $v$ is a solution of the equation $-\Delta v=\mu^2 v+f$. Thus the solution $u$ behaves asymptotically as the steady-state oscillation.

Turning to the equation (1) we have first to identify a counterpart of the oscillation $t\mapsto e^{-i\mu t}$. It is remarkable that an appropriate system of functions is already known. It has the form $\varphi_\omega (t)=E_\alpha (i^\alpha \omega^\alpha t^\alpha)$ where, as before, $1<\alpha <2$, $E_\alpha$ is the Mittag-Leffler function, $\omega>0$. For these functions, $\D \varphi_\omega =i^\alpha \omega^\alpha \varphi_\omega$. As $t\to \infty$, $\varphi_\omega (t)=\frac1\alpha e^{i\omega t}+O(t^{-\alpha })$ (see \cite{Dj,KST}; we follow the notations from \cite{KST} different from those in \cite{Dj} or \cite{PS}). Note that the ray $\arg z=\dfrac{\pi \alpha}2$ is the only one, on which the function $E_\alpha (z)$ has an oscillatory character; it has no zeroes on this ray (\cite{PS}, Theorem 4.2.1). Moreover, there exists an extensive theory of integral transforms based on the kernel $\varphi_\omega$. In various aspects, this theory is parallel to the standard Fourier analysis; see \cite{Dj,Gu,Mart}.

Therefore the principle of limiting amplitude for the equation (1) is formulated as follows. We consider a solution of the Cauchy problem for the equation (1) with $n\ge 3$, $f(t,x)=E_\alpha (i^\alpha \omega^\alpha t^\alpha)F(x)$. Under certain assumptions on $F$ and the initial functions we prove that
\begin{equation}
\frac{u(t,x)}{E_\alpha (i^\alpha \omega^\alpha t^\alpha)}\longrightarrow v(x),\quad t\to \infty ,
\end{equation}
for every $x\in \mathbb R^n$, where $v$ is a solution of the equation $\Delta v-i^\alpha \omega^\alpha v=-F$. A result of this kind is obtained also for an abstract equation $\D u+Au=f$ where $A$ is a non-negative selfadjoint operator on a Hilbert space.

The pointwise stabilization theorem for a solution $u(t,x)$ of the Cauchy problem with the initial condition $u(0,x)=u^0(x)$ for the heat equation $\dfrac{\partial u}{\partial t}=\Delta u$ is formulated as follows (see \cite{ER} and the survey papers \cite{DR,De}).

Let $u^0$ be a continuous bounded function. The solution $u(t,x)$ tends, for every $x\in \Rn$, to a constant $c$, as $t\to \infty$, if and only if, for for every $x_0\in \Rn$,
\begin{equation}
\lim\limits_{R\to \infty}\frac1{|K_R(x_0)|}\int\limits_{K_R(x_0)}u^0(x)\,dx=c
\end{equation}
where $K_R(x_0)$ is a ball of radius $R$ centered at $x_0$, $|K_R(x_0)|$ is its volume.

Note that the wave equation does not possess this property; only certain means of a solution stabilize \cite{GM}. Stabilization properties resembling the above one hold for hyperbolic equations with dissipative terms \cite{K86,K87}.

In this paper we prove a stabilization property of solutions of the diffusion-wave equation (1) similar to that of the heat equation. Thus, in this respect the equation (1) is closer to parabolic equations.

\bigskip
\section{Principle of limiting amplitude}

Let us consider the equation (1) with $n\ge 3$, $f(t,x)=E_\alpha (i^\alpha \omega^\alpha t^\alpha)F(x)$, $\omega >0$, and the initial conditions
\begin{equation}
u(0,x)=u^0(x),\quad \frac{\partial u(0,x)}{\partial t}=u^1(x).
\end{equation}
We assume that the functions $F,u^0,u^1$ are bounded; $F$ is locally H\"older continuous; $u^0\in C^1$ and its first derivatives are bounded and H\"older continuous with the exponent $\gamma >\dfrac{2-\alpha}\alpha$; $u^1$ is H\"older continuous. Under these assumptions, the problem (1),(4) possesses a classical solution $u(t,x)$ \cite{Ps,K14}. This means that $u(t,x)$ belongs to $C^2$ in $x$ for each $t>0$; $u(t,x)$ belongs to $C^1$ in $(t,x)$, and for any $x\in \Rn$ the fractional integral
$$
\left( I_{0+}^{2-\alpha }u\right) (t,x)=\frac{1}{\Gamma (2-\alpha
)}\int\limits_0^t(t-\tau )^{-\alpha +1}u'_\tau (\tau ,x)\,d\tau
$$
is continuously differentiable in $t$ for $t>0$;
$u(t,x)$ satisfies the equation and initial conditions.

Moreover, $u(t,x)$ admits the integral representation
\begin{multline}
u(t,x)=\int\limits_{\mathbb R^n} Z_1^{(\alpha )}(t,x;\xi )u^0(\xi )\,d\xi +\int\limits_{\mathbb R^n} Z_2^{(\alpha )}(t,x;\xi )u^1(\xi )\,d\xi \\
+\int\limits_0^td\tau \int\limits_{\mathbb R^n} Y(t-\tau,x;\xi)E_\alpha (i^\alpha \omega^\alpha \tau^\alpha)F(\xi )\,d\xi \overset{\text{def}}{=}u_1+u_2+u_3
\end{multline}
where the kernels satisfy the following estimates valid for $0<t<\infty$, $x\in \Rn$:
\begin{equation}
\left| Z_1^{(\alpha )}(t,x;\xi)\right| \le Ct^{-\alpha }|x-\xi |^{-n+2}\rho_\sigma (t,x,\xi );
\end{equation}
\begin{equation}
\left| Z_2^{(\alpha )}(t,x;\xi)\right| \le Ct^{-\alpha +1}|x-\xi |^{-n+2}\rho_\sigma (t,x,\xi );
\end{equation}
\begin{equation}
\left| Y^{(\alpha )}(t,x;\xi )\right| \le Ct^{\alpha -\frac{\alpha n}2-1}\mu_n(t^{-\alpha /2}|x-\xi |)\rho_\sigma (t,x,\xi ),
\end{equation}
where
$$
\rho_\sigma (t,x,\xi )=\exp \left\{ -\sigma (t^{-\alpha /2}|x-\xi |)^{\frac2{2-\alpha}}\right\},
$$
$$
\mu_n(z)=\begin{cases}
1, &\text{if $n=3$};\\
1+|\log z|, &\text{if $n=4$};\\
z^{-n+4}, &\text{if $n\ge 5$}.\end{cases}
$$
Here and below we denote by $C,\sigma$ various positive constants.

In fact, in the study of asymptotic properties of the function (5) it would be possible to remove the above smoothness assumptions considering $u$ as a kind of generalized solution.

\medskip
\begin{teo}
Assume, in addition to the boundedness, that $u^0,u^1,F\in L_1(\Rn )$, so that, in particular, $F\in L_2(\Rn )$. Then the limit relation (2) is valid for every $x\in \Rn$, with the function $v$ belonging to the Sobolev space $H^2(\Rn )$ and satisfying the equation
\begin{equation}
\Delta v-i^\alpha \omega^\alpha v=-F.
\end{equation}
\end{teo}

\medskip
{\it Proof}. It follows from (6) that
\begin{multline*}
|u_1(t,x)|\le Ct^{-\alpha }\int\limits_{\Rn}|x-\xi |^{-n+2}|u^0(\xi )\,d\xi \\
\le Ct^{-\alpha }\Biggl\{ \int\limits_{|x-\xi |\le 1}|x-\xi |^{-n+2}\,d\xi +\int\limits_{\Rn}|u^0(\xi )|\,d\xi \Biggr\} \le Ct^{-\alpha }\to 0,
\end{multline*}
as $t\to \infty$. Similarly, we find from (7) that $|u_2(t,x)|\le Ct^{-\alpha +1}\to 0$. Together with the asymptotics of the Mittag-Leffler function, this implies the relations
\begin{equation}
\frac{u_j(t,x)}{E_\alpha (i^\alpha \omega^\alpha t^\alpha)}\longrightarrow 0,\quad t\to \infty ,\quad j=1,2.
\end{equation}

While the qualitative behavior of the kernel $Y^{(\alpha )}$ is given by (8), to study $u_3$ we need an explicit representation \cite{Ps}:
\begin{multline*}
Y^{(\alpha )}(t-\tau ,x;\xi )=\Gamma_{\alpha ,n}(x-\xi ,t-\tau )\\
=2^{-n}\pi^{\frac{1-n}2}(t-\tau )^{\alpha -\frac{\alpha n}2-1}f_{\alpha /2}((t-\tau )^{-\alpha /2}|x-\xi |;n-1,\alpha -\frac{\alpha n}2).
\end{multline*}
Here
$$
f_{\alpha /2}(z;\mu ,\delta )=
\frac2{\Gamma (\mu /2)}\int\limits_1^\infty \Phi (-\alpha /2,\delta ,-zt)(t^2-1)^{\frac{\mu}2-1}dt, \quad \mu >0;
$$
$$
\Phi (-\alpha /2,\delta ;z)=\sum\limits_{m=0}^\infty \frac{z^m}{m!\Gamma (\delta -\frac{\alpha m}2)}
$$
is the Wright function.

It is known \cite{St} that
$$
\int\limits_0^\infty e^{-st}t^{\alpha -\frac{\alpha n}2-1}\Phi (-\alpha /2,\alpha -\frac{\alpha n}2;-\tau t^{-\alpha/2}|x|)\,dt=s^{-\alpha +\frac{\alpha n}2}e^{-\tau |x|s^{\alpha /2}}.
$$
Therefore we find that
$$
\int\limits_0^\infty e^{-st}\Gamma_{\alpha ,n}(x,t)\,dt
=\frac{2^{-n+1}\pi^{\frac{1-n}2}}{\Gamma (\frac{n-1}2)}s^{-\alpha +\frac{\alpha n}2}\int\limits_1^\infty (\tau^2-1)^{\frac{n-1}2-1}e^{-\tau |x|s^{\alpha /2}}d\tau ;
$$
the change in the order of integration was justified by an estimate of the Wright function given in Lemma 1 of \cite{Ps}.

Next we use the identity
$$
\int\limits_1^\infty (\tau^2-1)^{\nu -1}e^{-\mu \tau}d\tau =\frac1{\sqrt{\pi }}\left( \frac2{\mu}\right)^{\nu -\frac12}\Gamma (\nu )K_{\nu -\frac12}(\mu )
$$
$(\mu >0,\nu >0)$; see (\cite{GR}, 3.387.3). Here $K_{\nu -\frac12}$ is the Macdonald function. We come to the identity
\begin{equation}
\int\limits_0^\infty e^{-st}\Gamma_{\alpha ,n}(x,t)\,dt =2^{-n/2}\pi^{-n/2}|x|^{-\frac{n}2+1}s^{-\frac{\alpha}2 +\frac{\alpha n}4}K_{\frac{n}2-1}(s^{\alpha /2}|x|),\quad s>0,x\ne 0.
\end{equation}

For $x\ne 0$, the function $\Gamma_{\alpha ,n}(x,t)$ has no singularity in $t$ at the point $t=0$, due to the exponential decay of the function $\rho_\sigma$. As $t\to \infty$ and $x\ne 0$ fixed, we can estimate $\rho_\sigma$ from above by (1), and the inequality (8) provides the upper estimate of $\Gamma_{\alpha ,n}(x,t)$ by $\const \cdot t^{-1-\frac{\alpha}2}$, if $n=3$, $\const \cdot t^{-\alpha -1}|\log t^{-1}|$, if $n=4$, $\const \cdot t^{-\alpha -1}$, if $n\ge 5$. This means that the left-hand side of (11) is holomorphic in $s$ on the half-plane $\R s>0$ and continuous in $s$ on the closed half-plane $\R s\ge 0$. The function in the right-hand side of (11) possesses the same properties, so that (11) holds for $\R s\ge 0$. In particular, for any real $\omega$, we have the equality
\begin{equation}
\int\limits_0^\infty e^{-i\omega t}\Gamma_{\alpha ,n}(x,t)\,dt =2^{-n/2}\pi^{-n/2}|x|^{-\frac{n}2+1}(i\omega )^{-\frac{\alpha}2 +\frac{\alpha n}4}K_{\frac{n}2-1}((i\omega )^{\alpha /2}|x|).
\end{equation}

Let us write
$$
u_3(t,x)=\int\limits_\Rn k(t,x-\xi )F(\xi )\,d\xi
$$
where
$$
k(t,x)=\int\limits_0^t\Gamma_{\alpha ,n}(x,t-\tau )E_\alpha (i^\alpha \omega^\alpha \tau^\alpha)\,d\tau .
$$
Next we study the asymptotic behavior of $k(t,x)$, as $t\to \infty$. We begin with the estimate
\begin{multline*}
|k(t,x)|\le  C\int\limits_0^t (t-\tau )^{\alpha -\frac{\alpha n}2-1}\mu_n((t-\tau )^{-\alpha /2}|x|)e^{-\sigma ((t-\tau)^{-\alpha /2}|x|)^{\frac2{2-\alpha}}}d\tau \\
=C|x|^{-n+2}\int\limits_0^{t|x|^{-2/\alpha}}\theta^{\alpha -\frac{\alpha n}2-1}\mu_n(\theta^{-\alpha /2})e^{-\sigma \theta^{-\frac{\alpha}{2-\alpha}}}d\theta.
\end{multline*}
Estimating the last integral by the one over $(0,\infty )$ we find that
\begin{equation}
|k(t,x)|\le  C|x|^{-n+2}
\end{equation}
where $C$ does not depend on $t$.

Let us write the asymptotic formula for the Mittag-Leffler function \cite{Dj,KST} in the form
$$
E_\alpha (i^\alpha \omega^\alpha \tau^\alpha)=\frac1{\alpha}e^{i\omega \tau}+r(\tau )
$$
where $|r(\tau )|\le C(1+\tau )^{-1}$, $\tau \ge 0$. Thus
$$
k(t,x)=\frac1\alpha \int\limits_0^t\Gamma_{\alpha ,n}(x,t-\tau )e^{i\omega \tau}d\tau+\int\limits_0^t\Gamma_{\alpha ,n}(x,t-\tau )r(\tau )\,d\tau \overset{\text{def}}{=}k_1(t,x)+k_2(t,x).
$$

We have
$$
k_1(t,x)=\frac1\alpha e^{i\omega t}\int\limits_0^t\Gamma_{\alpha ,n}(x,\tau )e^{-i\omega \tau}d\tau .
$$
Since
$$
\frac{k_1(t,x)}{E_\alpha (i^\alpha \omega^\alpha t^\alpha)}\sim \frac{k_1(t,x)}{\frac1\alpha e^{i\omega t}},\quad t\to \infty ,
$$
we get the relation
\begin{equation}
\frac{k_1(t,x)}{E_\alpha (i^\alpha \omega^\alpha t^\alpha)}\longrightarrow \int\limits_0^\infty\Gamma_{\alpha ,n}(x,\tau )e^{-i\omega \tau}d\tau ,
\end{equation}
as $t\to \infty$, for any $x\ne 0$.

On the other hand,
\begin{multline*}
|k_2(t,x)|\le C\int\limits_0^t (t-\tau )^{\alpha -\frac{\alpha n}2-1}\mu_n((t-\tau )^{-\alpha /2}|x|)e^{-\sigma ((t-\tau)^{-\alpha /2}|x|)^{\frac2{2-\alpha}}}(1+\tau )^{-1}d\tau \\
\le C(1+\frac{t}2)^{-1} \int\limits_0^{t/2} \tau^{\alpha -\frac{\alpha n}2-1}\mu_n(\tau^{-\alpha /2}|x|)e^{-\sigma (\tau^{-\alpha /2}|x|)^{\frac2{2-\alpha}}}d\tau \\
+C\int\limits_{t/2}^t \tau^{\alpha -\frac{\alpha n}2-1}\mu_n(\tau^{-\alpha /2}|x|)e^{-\sigma (\tau^{-\alpha /2}|x|)^{\frac2{2-\alpha}}}d\tau \\
=C(1+\frac{t}2)^{-1}|x|^{-n+2}\int\limits_0^{\frac12 t|x|^{-2/\alpha}} \theta^{\alpha -\frac{\alpha n}2-1}\mu_n(\theta^{-\alpha /2})e^{-\sigma \theta^{-\frac{\alpha}{2-\alpha}}}d\theta \\
+C\int\limits_{\frac12 t|x|^{-2/\alpha}}^{t|x|^{-2/\alpha}} \theta^{\alpha -\frac{\alpha n}2-1}\mu_n(\theta^{-\alpha /2})e^{-\sigma \theta^{-\frac{\alpha}{2-\alpha}}}d\theta ,
\end{multline*}
so that
\begin{equation}
\frac{k_2(t,x)}{E_\alpha (i^\alpha \omega^\alpha t^\alpha)}\longrightarrow 0,\quad \text{as $t\to \infty$.}
\end{equation}

By (14) and (15),
$$
\frac{k(t,x)}{E_\alpha (i^\alpha \omega^\alpha t^\alpha)}\longrightarrow \int\limits_0^\infty\Gamma_{\alpha ,n}(x,\tau )e^{-i\omega \tau}d\tau ,\quad \text{as $t\to \infty$.}
$$
Now we recall (10) and use the dominated convergence theorem to show that
\begin{equation}
\frac{u(t,x)}{E_\alpha (i^\alpha \omega^\alpha t^\alpha)}\longrightarrow \int\limits_{\Rn} \left[ \int\limits_0^\infty\Gamma_{\alpha ,n}(x-\xi ,\tau )e^{-i\omega \tau}d\tau \right] F(\xi )\,d\xi .
\end{equation}
Our use of the dominated convergence theorem is based on the inequality (13) being performed in two stages. The domain of integration is decomposed into the union of $\{ \xi :\ |x-\xi |\le 1\}$ and $\{ \xi :\ |x-\xi |>1\}$. In the first integral we use the boundedness of $F$, while in the second integral we recall that $F\in L_1(\Rn )$.

By virtue of (12), the right-hand side of (16) equals $\int\limits_{\Rn}G(x-\xi )F(\xi )\,d\xi$ where
$$
G(z)=(2\pi)^{-n/2}|z|^{-\frac{n}2+1}(i\omega )^{-\frac{\alpha}2 +\frac{\alpha n}4}K_{\frac{n}2-1}((i\omega )^{\alpha /2}|z|).
$$
Comparing this with the expression for the Green function of the Laplacian on $\Rn$ (see \cite{T}, Section 13.7) we obtain the required limit relation (2) where the function $v$ satisfies (9).$\qquad \blacksquare$

\bigskip
\section{An equation with a selfadjoint operator}

In this section we consider the Cauchy problem
\begin{equation}
\left( \D u\right) (t)+Au(t)=f(t),\quad u(0)=u'(0)=0,
\end{equation}
where $A$ is a selfadjoint nonnegative operator (that is, $A\ge 0$) on a Hilbert space $\mathcal H$ with the inner product $(\cdot ,\cdot )$, $f(t)=E_\alpha (i^\alpha \omega^\alpha t^\alpha)f_0$, $f_0\in \mathcal H$. We will use the representation
$$
A=\int\limits_0^\infty \lambda \,d\mathcal E_\lambda
$$
where $\mathcal E_\lambda$ is the resolution of the identity (projection-valued spectral measure) corresponding to the operator $A$. Of course, our setting is a simple special case of the general theory of equations (17) (see \cite{LP}) guaranteeing, in particular, the uniqueness of a solution.

We call a function $u$ a strong solution of (17), if $u\in C^1([0,\infty),\mathcal H)$, the function $u$ takes values in $D(A)$, there exists a continuous derivative of the function $t\mapsto \int\limits_0^t(t-\tau )^{-\alpha +1}u'(\tau )\,d\tau$, and the equalities in (17) are satisfied.

\medskip
\begin{lem}
The strong solution of the problem (17) is given by the formula
\begin{equation}
u(t)=t^\alpha \int\limits_0^\infty \frac{\lambda E_{\alpha ,\alpha +1}(-\lambda t^\alpha)+i^\alpha \omega^\alpha E_{\alpha ,\alpha +1}(i^\alpha \omega^\alpha t^\alpha )}{\lambda +i^\alpha \omega^\alpha}\,d\mathcal E_\lambda f_0.
\end{equation}
Here $E_{\alpha ,\beta}(z)=\sum\limits_{k=0}^\infty \dfrac{z^k}{\Gamma (\alpha k+\beta )}$ is the Mittag-Leffler type function.
\end{lem}

\medskip
{\it Proof}. It is known (see Section 4.1.3 in \cite{KST}) that the function
$$
y_\lambda (t)=\int\limits_0^t (t-\tau )^{\alpha -1}E_{\alpha ,\alpha}(-\lambda (t-\tau ))\varphi (\tau )\,d\tau ,\quad \lambda >0,
$$
satisfies the scalar Cauchy problem
$$
\mathbb D^{(\alpha )}y_\lambda +\lambda y_\lambda =\varphi ,\quad y_\lambda (0)=y'_\lambda (0)=0.
$$
If $\varphi (\tau )=E_\alpha (i^\alpha \omega^\alpha \tau^\alpha)$, we use an integration formula from (\cite{Dj}, page 2) and find that
$$
y_\lambda (t)=\frac{\lambda E_{\alpha ,\alpha +1}(-\lambda t^\alpha)+i^\alpha \omega^\alpha E_{\alpha ,\alpha +1}(i^\alpha \omega^\alpha t^\alpha )}{\lambda +i^\alpha \omega^\alpha}t^\alpha .
$$

Therefore the function (18) satisfies (17), if it is possible to apply $\D$ under the sign of integral in (18). It is known (\cite{KST}, (1.8.27)) that $E_{\alpha ,\alpha +1}(z)=-\dfrac1z+O(|z|^{-2})$, $z\to -\infty$. This asymptotic relation shows that the function $(t,\lambda )\mapsto \lambda E_{\alpha ,\alpha +1}(-\lambda t^\alpha )$ is bounded on the set $\{ (t,\lambda ):\ |\lambda t^\alpha|>1\}$. On its complement, the function $E_{\alpha ,\alpha +1}(-\lambda t^\alpha)$ is obviously bounded. Therefore $|\lambda y_\lambda (t)|\le C$, so that the application of $\mathbb D^{(\alpha )}$ is legitimate. $\qquad \blacksquare$

For this setting, the principle of limiting amplitude is formulated as follows.

\medskip
\begin{teo}
Suppose that the point $\lambda =0$ is not an eigenvalue for the operator $A$, and that $f_0\in D(A^{-1})$. Let $u(t)$ be a strong solution of the problem (17) with $f(t)=E_\alpha (i^\alpha \omega^\alpha t^\alpha)f_0$. Then
\begin{equation}
\frac{u(t)}{E_\alpha (i^\alpha \omega^\alpha t^\alpha)}\longrightarrow (A+i^\alpha \omega^\alpha I)^{-1}f_0,\quad t\to \infty ,
\end{equation}
in the sense of convergence in $\mathcal H$.
\end{teo}

\medskip
{\it Proof}. Let us write
$$
u_1(t)=\int\limits_0^{t^{-\alpha}} \frac{\lambda E_{\alpha ,\alpha +1}(-\lambda t^\alpha)+i^\alpha \omega^\alpha E_{\alpha ,\alpha +1}(i^\alpha \omega^\alpha t^\alpha )}{\lambda +i^\alpha \omega^\alpha}t^\alpha \,d\mathcal E_\lambda f_0,
$$
$$
u_2(t)=\int\limits_{t^{-\alpha}}^\infty \frac{\lambda E_{\alpha ,\alpha +1}(-\lambda t^\alpha)+i^\alpha \omega^\alpha E_{\alpha ,\alpha +1}(i^\alpha \omega^\alpha t^\alpha )}{\lambda +i^\alpha \omega^\alpha}t^\alpha \,d\mathcal E_\lambda f_0.
$$

By the orthogonality property of spectral decompositions
$$
\|u_1(t)\|^2=\int\limits_0^{t^{-\alpha}} \left| \frac{\lambda E_{\alpha ,\alpha +1}(-\lambda t^\alpha)+i^\alpha \omega^\alpha E_{\alpha ,\alpha +1}(i^\alpha \omega^\alpha t^\alpha )}{\lambda +i^\alpha \omega^\alpha}t^\alpha \right|^2\,d(\mathcal E_\lambda f_0,f_0).
$$
The integrand is a bounded continuous function. Since the function $\lambda \mapsto (\mathcal E_\lambda f_0,f_0)$ is continuous at $\lambda =0$, it follows from the absolute continuity property of the Stieltjes integral (see Section 9.34 of \cite{Ka} or Section 20 of \cite{Sch}) that $u_1(t)\to 0$, as $t\to \infty$.

Next, for $|\lambda t^\alpha|>1$,
$$
\lambda E_{\alpha ,\alpha +1}(-\lambda t^\alpha)=t^{-\alpha}+r_1(t,\lambda )
$$
where $|r_1(t,\lambda )|\le C\lambda^{-1}t^{-2\alpha}$,
$$
i^\alpha \omega^\alpha E_{\alpha ,\alpha +1}(i^\alpha \omega^\alpha t^\alpha )=\frac1\alpha e^{i\omega t}-t^{-\alpha}+r_2(t)
$$
where $|r_2(t)|\le Ct^{-2\alpha}$ (see the asymptotics of $E_{\alpha ,\alpha +1}$ in \cite{KST} or \cite{Dj}). This results in the representation
$$
u_2(t)=\frac1\alpha e^{i\omega t}\int\limits_{t^{-\alpha}}^\infty \frac1{\lambda +i^\alpha \omega^\alpha}d\mathcal E_\lambda f_0+t^\alpha \int\limits_{t^{-\alpha}}^\infty \frac{r_1(t,\lambda )}{\lambda +i^\alpha \omega^\alpha}d\mathcal E_\lambda f_0
+t^\alpha r_2(t)\int\limits_{t^{-\alpha}}^\infty \frac1{\lambda +i^\alpha \omega^\alpha}d\mathcal E_\lambda f_0 .
$$

As $t\to \infty$,
$$
\frac1\alpha e^{i\omega t}\int\limits_{t^{-\alpha}}^\infty \frac1{\lambda +i^\alpha \omega^\alpha}d\mathcal E_\lambda f_0\longrightarrow \frac1\alpha e^{i\omega t}(A+i\alpha \omega^\alpha I)^{-1}f_0,
$$
\begin{multline*}
\left\| t^\alpha \int\limits_{t^{-\alpha}}^\infty \frac{r_1(t,\lambda )}{\lambda +i^\alpha \omega^\alpha}d\mathcal E_\lambda f_0 \right\|^2=\int\limits_{t^{-\alpha}}^\infty \left| \frac{t^\alpha r_1(t,\lambda )}{\lambda +i^\alpha \omega^\alpha}\right|^2 d(\mathcal E_\lambda f_0,f_0)\le Ct^{-2\alpha }\int\limits_0^\infty  \left| \frac{\lambda^{-1} }{\lambda +i^\alpha \omega^\alpha}\right|^2 d(\mathcal E_\lambda f_0,f_0)\\
=Ct^{-2\alpha }\left\| (A+i^\alpha \omega^\alpha)^{-1}A^{-1}f_0\right\|^2\to 0,
\end{multline*}
since $f_0\in D(A^{-1})$.

Finally,
$$
\left\| \int\limits_{t^{-\alpha}}^\infty \frac{t^\alpha r_2(t)}{\lambda +i^\alpha \omega^\alpha}d\mathcal E_\lambda f_0\right\|^2\le Ct^{-\alpha }\left\| (A+i^\alpha \omega^\alpha)^{-1}f_0\right\|^2\to 0,
$$
as $t\to \infty$. Using again the asymptotics of $E_{\alpha }(i^\alpha \omega^\alpha t^\alpha )$ we get the relation (19).$\qquad \blacksquare$

\medskip
For $\mathcal H=L_2(\Rn )$, $A=-\Delta$, with the domain $D(A)=H^2(\Rn )$ and $f_0\in H^2(\Rn )$, Theorem 2 gives the principle of limiting amplitude with the convergence in the sense of $L_2(\Rn )$.

\bigskip
\section{Stabilization}

Let $u_\alpha (t,x)$, $t\ge 0,x\in \Rn$, be a bounded solution of the Cauchy problem
\begin{equation}
\D u-\Delta u_\alpha =0,\quad u_\alpha (0,x) =u^0(x),\frac{\partial u_\alpha (0,x)}{\partial t}=0
\end{equation}
where $1<\alpha <2$, $u^0$ is a bounded continuous function. Under additional smoothness conditions, namely if $u^0$ is continuously differentiable, and its first derivatives are bounded and H\"older continuous with the exponent $\gamma >\dfrac{2-\alpha}\alpha$, the function
\begin{equation}
u_\alpha (t,x)=\int\limits_{\Rn}Z_1^{(\alpha )}(t,x;\xi )u^0(\xi )\,d\xi
\end{equation}
is the unique bounded classical solution of the problem (20); see \cite{K14,Ps}. Without the above smoothness assumptions, we can investigate the function (21) interpreting it as a generalized solution of the problem (20).

\medskip
\begin{teo}
The function $u_\alpha (t,x)$ possesses the property of pointwise stabilization: there exists a constant $c$, such that
\begin{equation}
u_\alpha (t,x)\longrightarrow c \quad \text{for any $x\in \Rn$, as $t\to \infty$},
\end{equation}
if and only if the initial function $u^0$ satisfies the condition (3).
\end{teo}

\medskip
{\it Proof}. Since
$$
\int\limits_{\Rn}Z_1^{(\alpha )}(t,x;\xi )\,d\xi =1
$$
(see \cite{Ps}), it is sufficient to consider the case where $c=0$.

Suppose that $u^0$ satisfies the condition (3) with $c=0$. Denote
$$
\left( V_xu^0\right) (r)=\frac1{|K_r(x)|}\int\limits_{K_r(x)}u^0(y)\,dy
$$
where $x\in \Rn$ is a fixed point. Given $\varepsilon >0$, there exists $N>0$, such that
\begin{equation}
\left| \left( V_xu^0\right) (r)\right| <\varepsilon \quad \text{for $r\ge N$}.
\end{equation}
Below we assume that $n\ge 3$. The cases $n=1$ and $n=2$ can be treated similarly.

The kernel $Z_1^{(\alpha )}$ has in fact the form \cite{Ps} $Z_1^{(\alpha )}(t,x;\xi )=H_\alpha (t,|x-\xi |)$ where
\begin{equation}
|H_\alpha (t,r)|\le Ct^{-\alpha }r^{-n+2}e^{-\sigma (t^{-\alpha /2}r)^{\frac2{2-\alpha}}};
\end{equation}
\begin{equation}
\left|\frac{\partial}{\partial r}H_\alpha (t,r)\right| \le Ct^{-\alpha }r^{-n+1}e^{-\sigma (t^{-\alpha /2}r)^{\frac2{2-\alpha}}}.
\end{equation}
In spherical coordinates,
$$
u_\alpha (t,x)=\int\limits_0^\infty H_\alpha (t,r)\,dr \int\limits_{S_r(x)}u^0(\omega )\,dS_r(\omega )
$$
where $S_r(x)$ is the sphere of radius $r$ centered at $x$. We write
$$
u_\alpha (t,x)=\int\limits_0^\infty H_\alpha (t,r)\left\{  \frac{\partial}{\partial r}\int\limits_0^r\left[ \int\limits_{S_\rho (x)}u^0(\omega )\,dS_\rho(\omega )\right] \,d\rho \right\} \,dr,
$$
integrate by parts and note that, due to (24) and the inequality
$$
\left| \int\limits_0^r\left[ \int\limits_{S_\rho (x)}u^0(\omega )\,dS_\rho(\omega )\right] \,d\rho \right| =\left| \int\limits_{K_r(x)}u^0(y)\,dy\right| \le Cr^n,
$$
the boundary terms equal zero. Thus,
$$
u_\alpha (t,x)=-\int\limits_0^\infty \frac{\partial H_\alpha (t,r)}{\partial r}\,dr \int\limits_{K_r(x)}u^0(y)\,dy=\const \cdot \int\limits_0^\infty r^n\frac{\partial H_\alpha (t,r)}{\partial r}\left(V_xu^0\right) (r)\,dr=I_1+I_2
$$
where $I_1$ and $I_2$ correspond to the integration over $(0,N)$ and $(N,\infty )$ respectively.

By (23) and (25),
$$
|I_2|\le Ct^{-\alpha}\varepsilon \int\limits_N^\infty re^{-\sigma (t^{-\alpha /2}r)^{\frac2{2-\alpha}}}dr=C\varepsilon \int\limits_{t^{-\alpha /2}N}^\infty se^{-\sigma s^{\frac2{2-\alpha}}}ds \le C\varepsilon \int\limits_0^\infty se^{-\sigma s^{\frac2{2-\alpha}}}ds,
$$
that is $|I_2|$ is small for all values of $t$. Meanwhile,
$$
|I_1|\le Ct^{-\alpha}\int\limits_0^N r\,dr <\varepsilon,
$$
as $t>t_0$, for some $t_0>0$. Therefore $u_\alpha (t,x)\to 0$, as $t\to \infty$, for each $x\in \Rn$.

Conversely, suppose that $u_\alpha (t,x)\to 0$, as $t\to \infty$, for every $x\in \Rn$. Let us consider, simultaneously with the problem (20), the Cauchy problem for the heat equation,
\begin{equation}
\frac{\partial u_1}{\partial t}-\Delta u_1=0,\quad u_1(0,x)=u^0(x),
\end{equation}
with the same initial function $u^0$. Let $Z^{(1)}(t,x-\xi )$ be the classical fundamental solution of this Cauchy problem.

It follows from the subordination identity \cite{Ba} that
\begin{equation}
Z^{(1)}(t,x-\xi )=t^{-\alpha }\int\limits_0^\infty \Phi_{1/\alpha }(st^{-\alpha })H_\alpha (s,|x-\xi |)\,ds
\end{equation}
where
$$
\Phi_\gamma (z)=\sum\limits_{m=0}^\infty \frac{(-z)^m}{m!\Gamma (-\gamma m+1-\gamma )},\quad 0<\gamma <1,
$$
is the Wright function. We have the inequality
\begin{equation}
\left| \Phi_{1/\alpha}(z)\right| \le Ce^{-\sigma z^\frac{\alpha}{\alpha -1}}
\end{equation}
(see the asymptotics of $\Phi_\gamma$ in \cite{KST}). Since $u^0$ is bounded, it follows from (27), (28) and estimates for fundamental solutions that
\begin{equation}
u_1(t,x)=t^{-\alpha }\int\limits_0^\infty \Phi_{1/\alpha }(st^{-\alpha })u_\alpha (s,x)\,ds
\end{equation}
for each $t>0$, $x\in \Rn$.

Let us transform (29) setting $t^\alpha =z$. Then
$$
u_1(z^{1/\alpha},x)=\frac1z \int\limits_0^\infty \Phi_{1/\alpha }\left( \frac{s}z\right) u_\alpha (s,x)\,ds=\int\limits_0^\infty \frac{s}z \Phi_{1/\alpha }\left( \frac{s}z\right) u_\alpha (s,x)\,\frac{ds}s.
$$

Recall that the expression
$$
\left( k\overset{\text{M}}{*}f\right) (\theta )=\int\limits_0^\infty k\left( \frac{\theta}\tau \right) f(\tau )\,\frac{d\tau}\tau
$$
is called the Mellin convolution (see, for example, \cite{BGT}); in fact, this is a convolution in the sense of harmonic analysis on the multiplicative group $(0,\infty )$. We have
$$
u_1(z^{1/\alpha},x)=\left( k\overset{\text{M}}{*}u_\alpha \right) (z,x)
$$
where $k(\zeta )=\zeta^{-1}\Phi_{1/\alpha }(\zeta^{-1})$.

The Mellin transform
$$
\check{k}(z)=\int\limits_0^\infty s^{z-1}k\left( \frac1s\right)\,ds=\int\limits_0^\infty s^z\Phi_{1/\alpha }(s)\,ds
$$
exists for $\R z>1$. By the Abelian theorem for Mellin convolutions (Theorem 4.1.6 in \cite{BGT}), the stabilization $u_\alpha (t,x)\to 0$, as $t\to \infty$, implies the relation $u_1(z^{1/\alpha},x)\to 0$, as $z\to \infty$, that is $u_1(t,x)\to 0$, as $t\to \infty$. This implies the limit relation (3) (with $c=0$), the necessary and sufficient condition of pointwise stabilization for the heat equation. $\qquad \blacksquare$

\section*{Acknowledgements}

This work was supported in part by Grant No. 03-01-12 of the National Academy of Sciences of Ukraine under the program of joint research projects with Siberian Branch, Russian Academy of Sciences.

The author is grateful to the anonymous referee for helpful comments.

\end{document}